\documentclass[11pt]{amsart}

\usepackage{amsmath}
\usepackage{amssymb}
\usepackage{graphicx}

\newtheorem{theorem}{Theorem}[section]

\newtheorem{lemma}[theorem]{Lemma}
\newtheorem{proposition}[theorem]{Proposition}
\theoremstyle{definition}
\newtheorem{definition}[theorem]{Definition}

\newtheorem{remark}[theorem]{Remark}

\def\be{\begin{equation}}
\def\ee{\end{equation}}
\def\bea{\be\begin{array}{rcl}}
\def\eea{\end{array}\ee}
\def\beas{\begin{eqnarray*}}
\def\eeas{\end{eqnarray*}}

\newcommand{\avg}[1]{\langle #1 \rangle}

\makeatletter
\def\proof{\@ifnextchar[\@xproof\@proof}
\def\@proof{{\noindent\it Proof.}\enspace}
\def\@xproof[#1]{{\noindent\it #1.}\enspace}
\makeatother

\def\th{\theta}

\newcommand{\E}{\operatorname{\mathbb{E}}}
\newcommand{\Cov}{\operatorname{Cov}}
\newcommand{\Var}{\operatorname{Var}}

\def\sqr#1#2{{\vbox{\hrule height.#2pt
     \hbox{\vrule width.#2pt height#1pt \kern#1pt
           \vrule width.#2pt}
     \hrule height.#2pt}}}
\def\square{\sqr74}
\def\proofend{{\unskip\nobreak\hfil\penalty50\hskip1em
             \hbox{}\nobreak\hfil\square \parfillskip=0pt
             \finalhyphendemerits=0 \par\goodbreak \vskip8mm}}

\numberwithin{equation}{section}

\title[Traveling waves in the McKean-Vlasov equation]{Traveling Waves in the McKean-Vlasov Equation under Sakaguchi-Kuramoto Interaction with Phase Frustration}

\author[J. Vukadinovic]{Jesenko Vukadinovic}

\address[J. Vukadinovic]{College of Staten Island CUNY, 2800 Victory Blvd, Staten Island, NY 10314, USA}
\email{{\tt jesenko.vukadinovic@csi.cuny.edu}}



\begin{document}
\begin{abstract}
We study the McKean–Vlasov equation for weakly coupled oscillators for the Sakaguchi–Kuramoto model. While the original Kuramoto model with purely sinusoidal coupling provides a good description for small densely connected networks, time delays in large networks generate symmetry-breaking phase offsets. Sakaguchi and Kuramoto proposed the simplest extension that captures this effect by incorporating a mean-field frustration parameter $\alpha$ into a single-mode interaction. We establish a continuous global phase transition from incoherence to a unique non-equilibrium traveling-wave state that takes the form of a rotating exponentially modified circular normal distribution. This is a novel skew extension of the von Mises distribution family -- parametrized by location parameter $\theta_0$, concentration parameter $r$, and skewness parameter $\beta$ -- that arises through exponential filtering of its Fourier spectrum. The extension is natural in that it preserves the fragile Bessel moment hierarchy, which ensures that the family remains globally identifiable -- a property not shared by existing skew extensions. The equation for traveling waves reduces to a mean-field self-consistency condition for $r$ and $\beta$. The latter plays a dual role, statistically as a skewness parameter and dynamically as the effective frustration in that $\tan\beta$ is the wave speed. Existence and uniqueness are proven by showing that the normalized mean resultant map -- the asymmetric deformation of the normalized Bessel ratio -- is strictly monotone along the isogones, rendering it a globally invertible map between natural (bare) and mean (effective) parameters. The proof combines tools from geometric function theory and analytic combinatorics.   
\end{abstract}

\maketitle

\section{Introduction}
Many natural oscillators are mathematically represented as globally attracting limit-cycles that encode the asymptotic behavior of nonlinear differential equations. When a large number of such oscillators interact, the resulting system of coupled subsystems is high-dimensional, rendering it analytically intractable. The phase reduction approach replaces each subsystem with one equation per oscillator that governs its phase. It was initially proposed by Winfree~\cite{Winfree}, who relied on separation between fast timescales for the relaxation to the limit-cycles and slow timescales for their mutual interaction. Kuramoto~\cite{Kuramoto} provided a rigorous justification by employing perturbative methods that apply when the coupling between the limit-cycle oscillators is weak enough relative to their intrinsic frequencies to allow for mutual phase perturbation without destabilization.

When unperturbed, each limit cycle is viewed as a uniform oscillator $\dot \theta_i=\omega_i$, 
where $\omega_i$ is the intrinsic frequency of the oscillation. Averaging methods yield the following system of $N$ coupled phase equations that govern the time evolution of $N$ interacting oscillators:
\be\label{DKuramoto}
\dot \theta_i=\omega_i+\frac\mu N\sum_{j=1}^N\Gamma_{ij}(\theta_j-\theta_i), \ \ \ i=1, \dots, N.
\ee
Each equation describes how a limit cycle adjusts its phase velocity according to input from other oscillators. The phase interaction functions (PIF), $\{\Gamma_{ij}\}$, are convolutions of perturbation functions and processes from the original limit cycle model~\cite{Ermentrout}. When the oscillators are identical and the interactions are global and equally weighted, then $\Gamma_{ij}=\Gamma=\partial_\theta W$, where $W$ is the so-called interaction potential. The parameter $\mu>0$ represents the interaction strength. 

With the aim of achieving maximal tractability while capturing the essence of spontaneous synchronization, Kuramoto proposed purely sinusoidal coupling between identical oscillators, $\Gamma(\theta)=\sin\theta$. It mimics how oscillators adjust their phase velocities to each other: it vanishes when oscillators are in phase or anti-phase, it pushes oscillators together when they are in near phase and apart when they are in near anti-phase. Furthermore, a single-mode interaction permits a mean-field reformulation, in which the microscopic interactions can be captured by macroscopic mean-field quantities -- the order parameter, and the average phase. 
Although highly symmetric and idealized, the Kuramoto model is considered the paradigm of synchronization. Its central feature is a phase transition from incoherence to coherence at a threshold $\mu_c$ that emerges as a result of the competition between the disorder encoded in the frequency distribution $g(\omega)$ and the order induced by the coupling. This is best captured by the following classical result: For $N\gg1$, for unimodal symmetric frequency distributions $g(\omega)$, it predicts a phase transition at $\mu_c=\frac2{\pi g(0)}$ from the incoherent state to a coherent steady state in which a single cluster of oscillators spontaneously synchronizes (see Strogatz et al. \cite{Strogatz, StrogatzMirollo}).

The Kuramoto model requires modifications to improve its physical and biological plausibility, as well as to give rise to diverse dynamical phenomena that go well beyond the stationary mono-stable single-cluster coherence. Neurobiological applications in particular are abundant with phase-reduced models. The visual cortex is an intricate network featuring dynamical phenomena that cannot be accounted for in the Kuramoto framework (e.g. multi-stable multi-cluster states and periodic switching between them, dynamical instabilities, traveling waves, disorder in spatial patterns). The main shortcoming of the Kuramoto model from the neurobiological perspective lies in its lack of explicit spatial embedding. Global and uniform coupling of the system is a good approximation in a small densely connected network, but it loses its validity at the scale of cortical circuits in large neural networks. A crucial step toward neurobiological plausibility is the incorporation of time delays in the PIF that originate from two sources: axonal conduction delays, arising from the finite propagation speed of action potentials along nerve fibers, and synaptic delays due to the time required for neurotransmitter release, diffusion across the synaptic cleft, and binding to the postsynaptic receptor. Axonal conduction delays in particular scale with the distance between neurons, leading to phase lags in spatially extended networks that accumulate and grow with network size. At the microscopic level, this cumulative effect manifests as distance-dependent phase offsets in the PIF, requiring explicit spatial embedding. At the macroscopic level, a mean-field approximation captures this effect through a single phase frustration parameter $\alpha$ in a one-mode interaction $\Gamma(\theta)=\sin(\theta-\alpha)$, while keeping the coupling global and uniform and avoiding the need for spatial embedding. It was proposed by Sakaguchi and Kuramoto~\cite{Sakaguchi1}, yielding the simplest extension of the classical Kuramoto model that is not a perturbation, but a physically generic situation in large networks.

Phase offsets have significant dynamical implications -- they introduce symmetry-breaking frustrations into the model, causing the interaction functions to push phases away from absolute synchrony and giving rise to traveling-wave dynamics. Other more intricate dynamics are also predicted, such as spatial wave patterns reminiscent of those observed in cortical activity~\cite{zanette}. While the model has attracted renewed interest in recent years (see~\cite{Acebron, brede} and references therein), a rigorous mathematical treatment of the phase transition from incoherence to a non-equilibrium traveling-wave regime is still lacking.  

The mean-field (continuum) limit of infinitely many weakly coupled noisy oscillators lends itself much more readily to rigorous mathematical analysis than the discrete model (\ref{DKuramoto}). The first modification is to encode the disorder — present in the discrete model through the heterogeneous distribution of intrinsic frequencies — by replacing it with Brownian noise added to the phase equations of identical oscillators. This serves a dual purpose: it reflects the intrinsic stochasticity of biological oscillators, and it introduces an entropy effect that makes the mean-field phase transition mathematically precise. The mean-field limit makes this exact — it replaces the high-dimensional stochastic particle system by a single deterministic PDE of Fokker–Planck type for the probability density, rendering the competition between entropy encoded in the diffusion term and coupling energy encoded in the drift term explicit, giving rise to the phase transition.

Sakaguchi and Kuramoto~\cite{Sakaguchi} extended the Kuramoto model to include stochastic fluctuations in the intrinsic frequencies. Adding a noisy forcing term to each deterministic equation in the phase-reduced model (assuming $\omega_i=0$), one arrives at its Langevin counterpart,
\be\label{sakaguchi}
dX_t^i=\frac \mu N\sum_{j=1}^N\Gamma(X_t^j-X_t^i)dt+\sqrt{2D}dB_t^i, \ \ \ i=1, \dots, N,
\ee
where $B_t^i$ are $N$ independent Brownian motions. Provided that $\Gamma$ is smooth, it is known~\cite{Oelschlager} that we can pass to the so-called mean-field (thermodynamic, continuum) limit. Consider the empirical measure
\[
\rho^{(N)}=\frac1N\sum_{i=1}^N\delta_{X_t^i}.
\]
Then $\mathbb{E}\left(\rho^{(N)}\right)$ converges as $N\to\infty$ in the sense of weak convergence of probability measures to some circular probability measure $\rho$ satisfying the following nonlinear nonlocal parabolic PDE:
\be\label{VMK}
\partial_t\rho =D\partial^2_\theta\rho+\mu\partial_\theta\left(\rho \partial_\theta (W\star \rho)\right).
\ee
Equation (\ref{VMK}) can also be understood as the Fokker-Planck equation for the McKean SDE
\[
\dot X_t=\mu \left(W'\star\rho\right)(X_t,t)+\sqrt{2D}dB_t,
\]
where $\rho=\rm{Law}(X_t)$. It also arises as the overdamped limit of the Vlasov-Fokker-Planck equation.  

In~\cite{ConstantinVukadinovic}, P.~Constantin and J.~V. considered the McKean-Vlasov equation with Kuramoto interaction, $W(\theta)=-\cos\theta$, and proved that at $\mu_c=2$, equation (\ref{VMK}) undergoes a continuous global phase transition. For $\mu\leq \mu_c$, the incoherent state is the unique solution, while for $\mu>\mu_c$, there is an additional, unique up to rotation, coherent equilibrium steady state that takes the form of a von Mises probability distribution. {\em `Continuous'} refers to the fact that the coherent state emanates at the bifurcation point from the incoherent state through a pitchfork bifurcation for the concentration parameter, and {\em `global'} to the fact that the solution persists for all $\mu>\mu_c$. The von Mises distribution family -- parametrized by location and concentration parameters -- is a regular and globally identifiable family by virtue of the fact that the {\em mean resultant map} -- the first-order modified Bessel ratio $\frac{I_1}{I_0}$ -- is injective. The equation for finding coherent equilibrium steady states reduces to a mean-field self-consistency condition for the family, which takes the form of the transcendental fixed point problem for the mean resultant map: $\mu \frac{I_1(r)}{I_0(r)}=r$.
Existence and uniqueness of the solution for $\mu>2$ are therefore proven by showing that the {\em normalized mean resultant map} -- the first order normalized modified Bessel ratio --  
\[
\mathbb{R}_{>0}\owns r\mapsto R=\frac{2I_1(r)}{rI_0(r)}\in(0,1)
\]
is a diffeomorphism.  

The paper at hand tackles the phase transition for the McKean-Vlasov equation with Sakaguchi-Kuramoto interaction $W(\theta)=-\cos(\theta-\alpha)$, $\alpha\in(-\pi/2,\pi/2)$. We extend the result of~\cite{ConstantinVukadinovic} to arrive at an analogous conclusion with the caveat that the equilibrium dynamics is replaced by a non-equilibrium one: at $\mu_c=2\sec\alpha$, equation (\ref{VMK}) undergoes a continuous, global phase transition from incoherence to a non-equilibrium steady state. For $\mu\leq \mu_c$, the incoherent state is the unique solution, while for $\mu>\mu_c$, there is an additional unique traveling-wave solution that takes the form of a rotating exponentially modified circular normal distribution. It is a novel skew extension of the von Mises distribution family -- parametrized by location parameter $\theta_0$, concentration parameter $r$, and skewness parameter $\beta$ -- that arises through multiplication by a skewness factor that is obtained by exponential filtering of its Fourier spectrum. The extension is natural in that it preserves the fragile Bessel moment hierarchy, rendering the family globally identifiable in parameters $r$ and $\beta$. We shall refer to it as the Exponentially Modified von Mises Distribution (EMvM), which, to the best of our knowledge, has not been previously introduced in the literature. The construction is the circular analogue of the exponentially modified Gaussian distribution, widely used in reaction time modeling and signal processing. The relationship between the EMvM distribution and the ex-Gaussian is canonical rather than merely analogical — in the large concentration limit the EMvM distribution converges to the ex-Gaussian, precisely as the von Mises distribution converges to the Gaussian.

The equation for finding traveling waves reduces to a mean-field self-consistency condition — a transcendental fixed point problem in $r$ and $\beta$, which plays a dual role, statistically as a skewness parameter and dynamically as the effective (observable) frustration in that $k=\tan\beta$ is the wave speed. Existence and uniqueness are proven by showing that the {\em normalized mean resultant map} -- the asymmetric deformation of the normalized Bessel ratio --
\[
 \mathbb{R}_{>0}\times\mathbb{R}\owns re^{i\beta}\mapsto\frac{2{I}_{1-i\tan\beta}(r)}{r{I}_{-i\tan\beta}(r)}=Re^{i\alpha}\in \mathcal{W}
\]
where 
\[
\mathcal{W}=\left\{z=Re^{i\alpha}\ \Big|\ \alpha\in(-\pi/2,\pi/2),\ 0<R<\cos\alpha\right\}=\left\{z:|z-\tfrac12|<\tfrac12\right\}
\]
is a disk, is a globally invertible map between natural (bare) and mean (effective) parameters. The proof of the global invertibility is the technical heart of the paper. The proof requires a change of variables to isogonal coordinates, reducing the problem to monotonicity statements for ratios of generalized (double) power series.  

Various works dealing with the phase transition for different interaction potentials exist already~\cite{CarrilloGvalaniPavliotisSchlichting, ChayesPanferov, ConstantinVukadinovic, Dawson, LuciaVukadinovic, Tamura}. They require the interaction potential to possess a symmetry — either oddness~\cite{Tamura} (considered nonphysical), or evenness (considered physical). Furthermore, the vast majority of them are based on tools from nonlinear functional analysis, in particular Crandall–Rabinowitz type arguments for the local analysis, and the so-called Rabinowitz alternative for the global analysis. The most comprehensive paper up to date is~\cite{CarrilloGvalaniPavliotisSchlichting}, which contains explicit criteria for the existence of both continuous and discontinuous phase transitions based on the Fourier coefficients of the interaction potential for quite general interactions by exploiting the symmetry in the problem. Applications include the Kuramoto model, the Hegselmann–Krausse model for opinion dynamics, the Keller–Segel model for bacterial chemotaxis, the Onsager model for liquid crystal alignment, and the Barré–Degond–Zatorska model for interacting dynamical networks.

The Sakaguchi–Kuramoto interaction $W(\theta)=-\cos(\theta-\alpha)$ is not even for $\alpha\not=0$, placing it outside the scope of these methods. The Ott–Antonsen reduction~\cite{OttAntonsen} provides a powerful formal tool for analyzing the collective dynamics of coupled oscillators, and has been applied broadly to identify traveling-wave solutions and their stability in Kuramoto-type models. However, the Ott–Antonsen approach yields formal rather than rigorous results and does not address the global structure of the bifurcation diagram. The present paper provides the first rigorous treatment of the phase transition for the asymmetric Sakaguchi–Kuramoto interaction, establishing global existence and uniqueness of the traveling-wave solution for all $\mu>\mu_c$.
\section{McKean-Vlasov equations for finding equilibrium and traveling-wave steady-states} 
Setting for convenience $D=1$ and introducing the mean-field potential, 
\[
V=\mu W\star\rho, 
\]
and the chemical potential, 
\[
\xi=\log\rho+V, 
\]
we rewrite~\eqref{VMK} as 
\be\label{ev}
\partial_t\rho=\partial_\theta\left(\rho\partial_\theta\xi\right). 
\ee
Since $\rho$ represents a probability distribution function on the unit circle $\mathbb{S}^1$, it is subject to positivity $\rho>0$, periodic boundary condition $\rho(0)=\rho(2\pi)$, and normalization condition $\int_{0}^{2\pi}\rho(\theta)\ d\theta=1$. Additionally, we introduce the dual density $\psi$ through the nonlinear nonlocal transformation:
\[
\rho\mapsto \psi=\rho\exp (\mu W\star\rho). 
\]
The equation for finding equilibrium steady-states then reads  
\[
\rho\partial_\theta\xi=c, 
\]
for some constant $c\in\mathbb{R}$. In terms of the dual density $\psi$, it becomes   
\[
\partial_\theta\psi=c\exp (\mu W\star\rho).  
\]
The periodic boundary condition implies $c=0$, and the normalization condition renders $\psi=Z^{-1}$, where $Z$ is the normalization constant, yielding the following McKean-Vlasov equation for finding equilibrium steady-states:  
\begin{align}\label{steadystates}
\rho=Z^{-1}\exp(-\mu W\star \rho),\\
Z=\int_{0}^{2\pi} \exp(-\mu W\star\rho)\ d\theta.\nonumber
\end{align}
In order to derive the equation for finding non-equilibrium steady-states in the form of traveling waves, we set out from the standard ansatz $\rho(t,\theta)=\rho(\theta-kt)$. Upon substituting it into (\ref{ev}), we obtain 
\[
\partial_\theta(k\rho+\rho\partial_\theta\xi)=0,
\]
so the McKean-Vlasov equation for finding traveling-wave steady-states becomes 
\be\label{trwaves}
\partial_\theta\rho+(\mu W'\star\rho+k)\rho=c,
\ee
 where $k\in\mathbb{R}$ and $c\in\mathbb{R}$ are constants to be determined from the periodic boundary condition $\rho(0)=\rho(2\pi)$ and the normalization condition. In terms of the dual density $\psi$, equation (\ref{trwaves}) reads 
\be\label{trawaves1}
\partial_\theta\psi+k\psi=c\exp (\mu W\star\rho). 
\ee

Equations (\ref{steadystates}) and (\ref{trwaves}) are nonlinear and nonlocal differential-integral equations. For interaction potentials with a finite number of Fourier modes, they permit a finite-dimensional reduction rendering them more tractable. Here, we consider the simplest case of a one-mode interaction potential with a phase shift, $W(\theta)=-\cos(\theta-\alpha)$, $\alpha\in(-\pi/2,\pi/2)$.   

For a given PDF $\rho$ on the unit circle $\mathbb{S}^1$, we introduce the notation for averages against that PDF: 
\[
\avg{f}_\rho=\int_{0}^{2\pi}f(\theta)\rho(\theta)\ d\theta.
\]
We introduce the order parameter, $r_0$, measuring the degree of coherence of $\rho$ and the average phase, $\theta_0\in[0,2\pi)$, via the equation 
\[
r_0 e^{\theta_0 i}=\avg{e^{\theta i}}_\rho. 
\]
Introducing the notation $r=\mu r_0$, we can rewrite it as 
\[
r =\mu \avg{e^{(\theta -\theta_0)i}}_\rho. 
\]
Equating the real and imaginary parts yields 
\begin{align}\label{OP}
r &=\mu \avg{\cos(\theta-\theta_0)}_\rho\\
0&=\avg{\sin(\theta-\theta_0)}_\rho. \nonumber
\end{align}
The mean-field potential $V=\mu W\star\rho$ can then be rewritten via the mean-field parameters $r$ and $\theta_0$ as 
\beas
V(\theta)=&=&-\mu\int_{0}^{2\pi}\cos(\theta-\theta_0+\theta_0-\tilde\theta-\alpha)\rho(\tilde\theta)\ d\tilde\theta\\
&=&-\mu[\cos(\theta-\theta_0-\alpha)\avg{\cos(\theta-\theta_0)}_\rho+\sin(\theta-\theta_0-\alpha)\avg{\sin(\theta-\theta_0)}_\rho]\\
&=&-r\cos(\theta-\theta_0-\alpha). 
\eeas
Introducing $\phi=\theta-\theta_0-\alpha$, (\ref{OP}) becomes 
\begin{align}\label{OP1}
r &=\mu\left(\cos\alpha \avg{\cos\phi}_\rho-\sin\alpha \avg{\sin\phi}_\rho\right)\\
0&=\cos\alpha \avg{\sin\phi}_\rho+\sin\alpha \avg{\cos\phi}_\rho. \nonumber
\end{align}
The stationary equation (\ref{steadystates}) then becomes 
\begin{align}\label{SMVE}
\rho(\theta)&=Z^{-1}\exp(r\cos\phi)\\
Z&=\int_{0}^{2\pi} \exp(r\cos\phi)\ d\phi, \nonumber 
\end{align}
where $r$ is to be determined from (\ref{OP1}).  

The equation for traveling waves (\ref{trwaves}), on the other hand, becomes
\be\label{trwaves2}
\partial_\phi\rho+r\sin\phi \rho+k\rho=c, 
\ee
where $r$, $k$ and $c$ are to be determined from (\ref{OP1}) and periodic and normalization conditions. Rewritten in terms of the dual density $\psi$, it becomes 
\be\label{trawaves3}
\partial_\phi\psi+k\psi=c\exp(-r\cos\phi). 
\ee

It turns out that (\ref{OP1}) can be transformed into two transcendental equations in unknowns $r$ and $k$, involving appropriately asymmetrically extended modified Bessel functions. The PDF $\rho$ is then found in the form of an appropriately asymmetrically extended von Mises PDF with $r$ representing the concentration, and $k$ a measure of skewness. These extensions will be introduced in the next sections. The constant $c$ is subsequently determined from the normalization via the formula 
\be\label{c}
c=\frac1{2\pi}(r\avg{\sin\phi}_\rho+k).
\ee
\section{Von Mises distribution and modified Bessel functions}
We begin by recalling the following: 
\begin{definition}
The von Mises distribution family (vM) on the circle $\mathbb{S}^1\cong[0,2\pi)$ is the two-parameter family of probability density functions:
\[
\rho_{vM}(\theta;\theta_0,r)=\frac{e^{r\cos(\theta-\theta_0)}}{2\pi I_0(r)},\ \theta\in[0,2\pi), 
\]
parametrized by the location parameter $\theta_0\in [0,2\pi)$, and the concentration parameter $r>0$, where 
\be\label{MBF}
I_n(r)=\frac1{2\pi}\int_{0}^{2\pi}\cos n\phi \exp(r\cos\phi)\ d\phi, \ n\in\mathbb{Z}, \ \ r\in\mathbb{R}
\ee
are the modified Bessel functions. In view of the Jacobi-Anger expansion for the unnormalized kernel of the distribution,  
\[
e^{r\cos\phi}=\sum_{n\in\mathbb{Z}}I_n(r)\cos n\phi,\ r,\phi\in\mathbb{R}, 
\]
the characteristic function of the von Mises distribution is given by the circular moments, 
\[
\varphi_{vM}(n;\theta_0,r)=\frac{I_n(r)}{I_0(r)}e^{in\theta_0}, \ n\in\mathbb{Z},\ \theta_0\in[0,2\pi),\ r>0. 
\]
\end{definition}
In the following proposition, we summarize some important properties of the modified Bessel functions of the first kind $I_\nu$, including their analytic continuation in the order $\nu$. (See~\cite{Amos, baricz, LuciaVukadinovic, Watson}, and references therein). 
\begin{proposition}
We have: 
\begin{enumerate}
\item For $\nu\in\mathbb{C}$, the modified Bessel Bessel function of the first kind, $I_\nu(z)$, is a particular solution of the second-order modified Bessel differential equation of the first kind, 
    \be\label{besseleqn}
    z^2y'' + zy' -(z^2 + \nu^2)y = 0, \ \nu, z\in\mathbb{C}.  
    \ee
    For $n\in\mathbb{Z}$, function $I_n $ has the integral representation (\ref{MBF}). 
    \item For $\nu\in\mathbb{C}$ and $z\in\mathbb{C}$, $I_\nu$ has the ascending power series representation
    \[
    I_\nu(z)=\sum_{p=0}^\infty\frac{1}{p!\,\Gamma(\nu+p+1)}\left(\frac{z}{2}\right)^{\nu+2p}, 
    \]
    which is entire in the order $\nu$ and reduces, for $\nu=n=0,1,2,\dots$, to
    \[
    I_n(z)=\sum_{p=0}^\infty\frac{1}{p!(n+p)!}\left(\frac{z}{2}\right)^{n+2p}. 
    \]
    \item For $\mu,\nu\in\mathbb{C}$ and $z\in\mathbb{C}$, the product $I_\mu I_\nu$ has the power series representation
    \be\label{prodseries}
    I_\mu(z)\,I_\nu(z)=\sum_{p=0}^\infty
    \frac{\Gamma(\mu+\nu+2p+1)}{p!\,\Gamma(\mu+p+1)\,\Gamma(\nu+p+1)\,\Gamma(\mu+\nu+p+1)}
    \left(\frac{z}{2}\right)^{\mu+\nu+2p}. 
    \ee
    In particular, its square has the power series representation 
    \be\label{quadratic}
    I_n(z)^2=\sum_{p=0}^\infty \frac{1}{(n+p)!^2}{{2n+2p}\choose{2n+p}}\left(\frac{z}{2}\right)^{2n+2p},\ \ n=0,1,2,\dots
    \ee
    and 
    \be\label{quadratic1}
    I_0(z)^2=\sum_{p=0}^\infty \frac {(2p)!}{p!^4}\left(\frac{z}{2}\right)^{2p},   
    \ee
    while the imaginary-order square is, for $k\in\mathbb{R}$,
    \be\label{imagsquare}
    |I_{-ik}(z)|^2=I_{-ik}(z)\,I_{ik}(z)=\sum_{p=0}^\infty\frac{(2p)!}{p!^2\,|\Gamma(p+1+ik)|^2}\left(\frac{z}{2}\right)^{2p}. 
    \ee
    \item The following symmetry relations hold for $n\in\mathbb{Z}$ and $z\in\mathbb{C}$: 
    \[
    I_{-n}(z)=I_n(z)
    \]
    and 
    \[
    I_n(-z)=(-1)^nI_n(z).
    \]
    \item The following two recursion relations hold for $\nu\in\mathbb{C}$, $z\in\mathbb{C}$:
    \be\label{rec1}
    2I_\nu'(z)=I_{\nu-1}(z)+I_{\nu+1}(z), 
    \ee
    and 
    \be\label{rec2}
    \frac{2\nu I_\nu(z)}{z}=I_{\nu-1}(z)-I_{\nu+1}(z). 
    \ee
    \item The function $g(z,t)=\exp\left(\frac z2(t+t^{-1})\right)$ is the generating function for modified Bessel functions of the first kind. More specifically, for $z\in\mathbb{C}$ and $t\in\mathbb{C}\backslash\{0\}$, 
    \[
    g(z,t)=\sum_{n\in\mathbb{Z}}I_n(z)t^n. 
    \]
    \item We have 
    \be\label{one}
    \sum_{n\in \mathbb{Z}}(-1)^n I_n(z)^2=1,
    \ee
    and 
    \be\label{two}
    \sum_{n\in \mathbb{Z}} I_n(z)^2=I_0(2z). 
    \ee
    \item For $\mu,\nu\in\mathbb{C}$ with $\operatorname{Re}(\mu+\nu)>-1$, the product $I_\mu I_\nu$ has the integral representation
    \be\label{prodintegral}
    I_\mu(z)\,I_\nu(z)=\frac{2}{\pi}\int_0^{\pi/2} I_{\mu+\nu}\big(2z\cos\varphi\big)\,\cos\big((\mu-\nu)\varphi\big)\,d\varphi. 
    \ee
    In particular, for $k\in\mathbb{R}$,
    \be\label{imagintegral}
    |I_{-ik}(z)|^2=\frac{2}{\pi}\int_0^{\pi/2} I_0\big(2z\cos\varphi\big)\,\cosh(2k\varphi)\,d\varphi. 
    \ee
\item The $n\textsuperscript{th}$ order modified Bessel ratio, 
\[
\mathbb{R}_{>0}\owns r\mapsto \frac{I_{n}(r)}{I_{n-1}(r)}\in(0,1),\ \ n=1,2,\dots
\]
is strictly increasing bijection and decreasing in $n$, while the $n\textsuperscript{th}$ normalized modified Bessel ratio
\[
\mathbb{R}_{>0}\owns r\mapsto \frac{2nI_{n}(r)}{rI_{n-1}(r)}\in(0,1)\ \ n=1,2,\dots
\]
is strictly decreasing bijection and increasing in $n$.   
\end{enumerate}
\end{proposition}
\begin{definition}
We refer to the modulus of the first circular moment (the first order Bessel ratio)
\[
\mathbb{M}:\mathbb{R}_{>0}\owns r\mapsto \frac{I_1(r)}{I_0(r)}\in(0,1)
\]
and the first order normalized Bessel ratio
\[
\mathcal{M}:\mathbb{R}_{>0}\owns r\mapsto R=\frac{2I_1(r)}{rI_0(r)}\in(0,1)
\]
as the mean resultant length, and the normalized mean resultant length, respectively.  
\end{definition}
\begin{remark}
The injectivity of the mean resultant map means that vM is globally identifiable and possesses a unique maximum likelihood estimate (MLE) for any non-uniform sample. The injectivity of the normalized mean resultant map ensures that the mean-field self-consistency equation has a unique solution, establishing a one-to-one correspondence between the natural parameter -- the interaction strength $\mu$ and the mean-field parameter -- the concentration parameter $r$.
\end{remark} 

Substituting (\ref{SMVE}) into (\ref{OP1}), the latter becomes the following system in terms of modified Bessel functions: 
\begin{align}\label{OP2}
r &=\mu \cos\alpha \frac{I_1(r)}{I_0(r)}\\
0&=\sin\alpha I_1(r). \nonumber 
\end{align}
This system is overdetermined: the second equation is only satisfied if $\alpha=0$, or $r=0$ (yielding the trivial solution $\rho=1/(2\pi)$). Therefore, for $\alpha\not=0$ the only possible stationary solution is the trivial one, while for $\alpha=0$, nontrivial solutions are obtained by solving the following mean-field self-consistency condition:  
\be\label{BR}
\mu \mathbb{M}(r)=r.  
\ee
Equation (\ref{BR}) undergoes a global supercritical pitchfork bifurcation at $\mu_c=2$: in addition to the trivial solution $r=0$, for $\mu>\mu_c$, it has a unique nontrivial solution $r=\mathcal{M}^{-1}(2/\mu)$. The coherent steady states are then given by $\rho(\theta)=\rho_{vM}(\theta;\theta_0,r)$ for any $\theta_0\in[0,2\pi)$, and they are unique up to the choice of $\theta_0$.

\section{Exponentially modified (skew) von Mises distribution}
Before we introduce the exponentially modified von Mises distribution, we first recall the wrapped exponential distribution. 
\begin{definition}
The wrapped exponential distribution family (WE) on the circle $\mathbb{S}^1\cong[0,2\pi)$ is the two-parameter family of probability density functions: 
\[
\rho_{WE}(\theta;\theta_0,k)=\frac{ke^{-k(\theta-\theta_0)}}{1-e^{-2k\pi}},\ \theta\in[0,2\pi), 
\]
parametrized by the location parameter $\theta_0\in[0,2\pi)$ and a decay rate parameter $k>0$. 
In view of the Fourier series expansion, 
\[
\frac1{2\pi}\sum_{n\in\mathbb{Z}}\frac{1}{1-in/k}e^{-in\phi}=\frac{ke^{-k\phi}}{1-e^{-2k\pi}},\ \phi\in[0,2\pi), 
\]
its characteristic function is given by 
\[
\varphi_{WE}(n;\theta_0,k)=\frac{e^{in\theta_0}}{1-\frac{ni}{k}}. 
\]
\end{definition}
For a given concentration parameter $r>0$ and a given decay rate parameter $k>0$, we now seek to obtain the unnormalized kernel $\chi_{k,r}$ of the exponentially modified von Mises PDF as the (unique) periodic solution of the equation 
\be\label{AEvMD}
\frac{d\chi}{d\phi}+(r\sin\phi +k)\chi=k. 
\ee
Equivalently, its skewing factor, $\psi_{k,r}=\chi_{k,r}e^{-r\cos\phi}$, is the (unique) periodic solution of 
\be\label{trawaves4}
\frac{d\psi}{d\phi}+k\psi=ke^{-r\cos\phi}. 
\ee
The Fourier coefficients of $\psi_{k,r}$,  
\[
a_n(k,r)=\frac1{2\pi}\int_{0}^{2\pi}\psi_{k,r}(\phi)\cos n\phi\ d\phi, \ \ b_n(k,r)=\frac1{2\pi}\int_{0}^{2\pi}\psi_{k,r}(\phi)\sin n\phi\ d\phi,\ n\in\mathbb{Z},
\]
are easily computed as 
\[ 
a_n(k,r)=\frac{k^2}{k^2+n^2}I_n(-r), \ \ b_n(k,r)=\frac{kn}{k^2+n^2}I_n(-r),\ n\in\mathbb{Z}. 
\]
Therefore, for the unnormalized kernel of the extension we have, 
\be\label{chi}
\chi_{k,r}(\phi)=e^{r\cos\phi}\psi_{k,r}(\phi), 
\ee
with the skewing factor, 
\begin{align}\label{psi}
\psi_{k,r}(\phi)&=\sum_{n\in\mathbb{Z}}I_n(-r)\frac{1}{1-in/k}e^{-in\phi}\\
&=\sum_{n\in\mathbb{Z}} I_n(-r)\Big[\frac{k^2}{k^2+n^2}\cos n\phi+\frac{nk}{k^2+n^2}\sin n\phi\Big]\nonumber. 
\end{align}
The product in the frequency space yields the following convolution representation of the solution in the phase space: \[
\psi_{k,r}(\phi)=\frac{k}{1-e^{-2k\pi}}\int_0^{2\pi}e^{-r\cos(\phi-\varphi)-k\varphi}\ d\varphi, 
\]
and 
\[
\chi_{k,r}(\phi)=\frac{k}{1-e^{-2k\pi}}\int_0^{2\pi}e^{r(\cos\phi-\cos(\phi-\varphi))-k\varphi}\ d\varphi.  
\]
For the normalization, we then have 
\beas
Z_{k,r}&=&\frac{k}{1-e^{-2k\pi}}\int_0^{2\pi}\int_0^{2\pi}e^{r(\cos\phi-\cos(\phi-\varphi))-k\varphi}\ d\varphi\ d\phi\\
&=&\frac{2\pi k}{1-e^{-2k\pi}}\int_0^{2\pi} I_0\big(2r\sin(\varphi/2)\big)e^{-k\varphi}\,d\varphi\\
&=&\frac{4k\pi}{\sinh k\pi}\int_0^{\pi/2} I_0\big(2r\cos\varphi\big)\cosh(2k\varphi)\,d\varphi\\
&=&\frac{2\pi^2 k}{\sinh k\pi}\,|I_{ik}(r)|^2\\
&=&2\pi|\Gamma(1+ik)|^2|I_{ik}(r)|^2.
\eeas 
[We first used the elementary identity $\int_0^{2\pi}e^{a\cos\phi+b\sin\phi}\,d\phi=2\pi I_0(\sqrt{a^2+b^2})$ applied with
$a=r(1-\cos\varphi)$, $b=-r\sin\varphi$, and $\sqrt{a^2+b^2}=2r\sin(\varphi/2)$, followed by the product integral~\eqref{imagintegral}, and the identity $|\Gamma(1+ik)|^2=\pi k/\sinh k\pi$.]
\begin{definition}
The exponentially modified von Mises distribution family (EMvM) on the circle $\mathbb{S}^1\cong[0,2\pi)$ is the two-parameter family of probability density functions:
\[
\rho_{EMvM}(\theta;\theta_0,k,r)=\frac{\chi_{k,r}(\theta-\theta_0)}{Z_{k,r}}
\]
parametrized by the location parameter $\theta_0\in[0,2\pi)$, the concentration parameter $r>0$, and the decay rate parameter $k>0$.  
\end{definition}
\begin{remark}
To the best of our knowledge, EMvM has not yet been formally introduced in the literature. It is the circular analogue of the Exponentially Modified Gaussian (EMG) or ex-Gaussian. In the context of circular statistics, $k$ provides a measure of skewness. We further introduce $\beta=\arctan (k)$, which provides a normalized measure of the distribution's deviation from the symmetric von Mises distribution. We also write $\rho_{EMvM}(\theta;\theta_0,\beta,r)$ as shorthand for $\rho_{EMvM}(\theta;\theta_0,\tan\beta,r)$. 
\end{remark}
\begin{remark}
    The unnormalized kernel $\chi_{k,r}$ and EMvM can be extended to $r\leq0$ and $k<0$ in a straightforward fashion. For $k=0$, equation (\ref{trawaves4}) no longer has a unique periodic solution, but rather it is satisfied by any constant. We can continuously extend $\psi_{k,r}$ by setting $\psi_{0,r}(\phi)=I_0(r)$ and $\chi_{0,r}(\phi)=I_0(r)e^{r\cos\phi}$. Consequently, $\rho_{EMvM}(\theta;\theta_0,0,r)=\rho_{vM}(\theta;\theta_0,r)$, so in this sense EMvM is a continuous extension of vM.  
\end{remark}  
\begin{definition}
For given $k,r\in\mathbb{R}$, $n\in \mathbb{Z}$, let 
\[
\mathbb{I}_n(k,r)=C_n(k,r)+iS_n(k,r)=\frac1{2\pi}\int_{0}^{2\pi}\chi_{k,r}(\phi)e^{-in\phi}\ d\phi
\]
be the moment coefficients of EMvM. 
\end{definition}
\begin{proposition}\label{prop:cont}
For $n\in\mathbb{Z}$ and $k,r\in\mathbb{R}$, the moment coefficients of EMvM have the closed form
\[
   \mathbb{I}_n(k,r)=|\Gamma(1+ik)|^2\,I_{ik}(r)\,I_{n-ik}(r),
\]
where $I_\nu(z)$ is the modified Bessel function of the first kind, analytically continued in the order $\nu$. In particular, the zeroth moment is
\[
   \mathbb{I}_0(k,r)=C_0(k,r)=|\Gamma(1+ik)|^2\,|I_{ik}(r)|^2,
\]
and the characteristic function of EMvM is given by the circular moments 
\[
\varphi_{EMvM}(n;\theta_0,k,r)=\frac{\mathbb{I}_n(k,r)}{\mathbb{I}_0(k,r)}e^{in\theta_0}=\frac{{I}_{n-ik}(r)}{{I}_{-ik}(r)}e^{in\theta_0}, \ n\in\mathbb{Z},\ \theta_0\in[0,2\pi), k,r\in\mathbb{R}. 
\]
\end{proposition}
\proof
The moment recursions arising from the kernel~\eqref{AEvMD} combine, for
$n\ge1$, into the single complex recursive relation
\[
\mathbb{I}_{n-1}(k,r)-\mathbb{I}_{n+1}(k,r)=\frac{2(n-ik)}{r}\,\mathbb{I}_n(k,r),
\]
which is exactly the modified-Bessel recurrence~\eqref{rec2} at $\nu=n-ik$, $z=r$. Both
$\{\mathbb{I}_n(k,r)\}$ and $\{I_{n-ik}(r)\}$ are recessive as $n\to+\infty$
($\mathbb{I}_n(k,r)\to0$ as Fourier coefficients of the integrable kernel $\chi_{k,r}$;
$I_{n-ik}(r)\sim(r/2)^{n-ik}/\Gamma(n-ik+1)\to0$). The recessive solution of a
second-order recursion is unique up to a constant, so $\mathbb{I}_n(k,r)=c(k,r)\,I_{n-ik}(r)$ with
$c(k,r)$ independent of $n$. Therefore 
\[
c(k,r)=\mathbb{I}_0(k,r)/I_{-ik}(r)=Z_{k,r}/(2\pi I_{-ik})=|\Gamma(1+ik)|^2 I_{ik}(r), 
\]
and the conclusion follows. 
\proofend

\begin{definition}
The first-order Bessel ratio and the first-order normalized Bessel ratio are respectively referred to as the mean resultant map and the normalized mean resultant map for EMvM. They are given by
\[
\mathbb{R}\times\mathbb{R}_+\owns(k,r)\mapsto\mathbb{M}(k,r)=\frac{\mathbb{I}_1(k,r)}{\mathbb{I}_0(k,r)}=\frac{{I}_{1-ik}(r)}{{I}_{-ik}(r)}
\]
and 
\[
\mathbb{R}\times\mathbb{R}_+\owns(k,r)\mapsto\mathcal{M}(k,r)=\frac{2\mathbb{I}_1(k,r)}{r\mathbb{I}_0(k,r)}=\frac{{2I}_{1-ik}(r)}{{rI}_{-ik}(r)}. 
\]
We further introduce the notation 
\[
T_1(k,r)=\frac{S_1(k,r)}{C_1(k,r)}.
\] 
We write 
\[
\mathcal{M}(k,r)=\mathcal{C}(k,r)+i\mathcal{S}(k,r)=\mathcal{R}(k,r)e^{i\mathcal{A}(k,r)}.
\]
The normalized mean resultant angle and length are, respectively, 
\begin{align*}
\mathcal{A}(k,r)&=\mathrm{Arg}(\mathcal{M}(k,r))=\arctan\left(T_1(k,r)\right),\\
\mathcal{R}(k,r)&=|\mathcal{M}(k,r)|=\sqrt{\mathcal{C}(k,r)^2+\mathcal{S}(k,r)^2}.  
\end{align*}
\end{definition}
\section{Statements of the main results}
\begin{remark}
Integration in~\eqref{AEvMD} yields
\[
r\int_0^{2\pi}\sin\phi \chi_{k,r}(\phi)\ d\phi+kZ_{k,r}=2\pi k, 
\]
so~\eqref{AEvMD} becomes 
\[
\frac{d\chi_{k,r}}{d\phi}+(r\sin\phi +k)\chi_{k,r}=\frac1{2\pi}\left(r\int_0^{2\pi}\sin\phi \chi_{k,r}(\phi)\ d\phi+kZ_{k,r}\right).  
\]
The density $\rho(\phi)=\frac{\chi_{k,r}(\phi)}{Z_{k,r}}$ then satisfies 
\[
\frac{d\rho}{d\phi}+(r\sin\phi +k)\rho=\frac1{2\pi}\left(r\int_0^{2\pi}\sin\phi \rho(\phi)\ d\phi+k\right), 
\]
which is exactly (\ref{trwaves2}) with $c$ given through formula~\eqref{c}.   

In summary, when $\alpha\not=0$ and $r\not=0$, a distribution $\rho$ satisfies (\ref{trwaves2}) if and only if $\rho(\theta)=\rho_{EMvM}(\theta;\theta_0,k,r)$, for some arbitrary $\theta_0\in[0,2\pi)$, and for $k\in\mathbb{R}$ and $r>0$ that satisfy the following mean-field self-consistency condition: 
\be\label{OP3}
\mu \mathbb{M}(k,r)=e^{i\alpha}r. 
\ee
\end{remark}
We now state two main theorems of the paper. The first theorem concerns the bifurcation diagram for the mean-field self-consistency condition (\ref{OP3}). In view of the remark above, the second theorem regarding the phase transition from incoherence to a unique traveling wave is an immediate consequence thereof.    
\begin{theorem}\label{main0}
For $\alpha\in(0,\pi/2)$ ($\alpha\in(-\pi/2,0)$), consider the mean-field self-consistency condition (\ref{OP3}). It undergoes a global supercritical pitchfork bifurcation in $\mu$ at the bifurcation point $\mu_\alpha=2\sec\alpha$. For $\mu\leq \mu_\alpha$, the trivial solution $r=0$ is the unique solution, while for $\mu>\mu_\alpha$, an additional solution $(\mu, k_\alpha(\mu),r_\alpha(\mu))$ emanates from the bifurcation point $(\mu_\alpha,k_\alpha(\mu_\alpha), r_\alpha(\mu_\alpha))=(2\sec\alpha,\tan\alpha,0)$. This bifurcation diagram is global in the sense that the curve consisting of nontrivial solutions, 
\[
\Gamma_{\alpha}=\{(k_\alpha(\mu),r_\alpha(\mu))\ \big|\ \mu> \mu_\alpha\},
\]
can also be parametrized by means of an increasing bijection $K_{\alpha}:(0,\infty)\to(\tan\alpha,\infty)$ (decreasing bijection $K_{\alpha}:(0,\infty)\to(-\infty,\tan\alpha)$):
\[
\Gamma_{\alpha}=\{(K_{\alpha}(r),r)\ \big|\ r>0\}, 
\]
exhausting all possible nontrivial solutions.  

Furthermore, for $\alpha\in(0,\pi/2)$,
\begin{align*}
&\lim_{\mu\to \mu_\alpha^+}k_\alpha(\mu)=\tan\alpha\\
&\lim_{\alpha\to\pi/2^-}k_\alpha(\mu)=\infty\\
&\lim_{\mu\to \infty}k_\alpha(\mu)=\infty\\
&\lim_{\alpha\to\pi/2^-}r_\alpha(\mu)=\infty \\
&\lim_{\mu\to \infty}r_\alpha(\mu)=\infty. 
\end{align*}
For $\alpha=\pm\pi/2$, the trivial solution is unique for all $\mu>0$.    
\end{theorem}
\begin{theorem}\label{main}
Consider equation~\eqref{trwaves} with a one-mode interaction potential:
\be\label{potential}
W(\theta)=-\cos(\theta-\alpha),\ \alpha\in(-\pi/2,\pi/2)
\ee
Equation~\eqref{trwaves} undergoes a continuous global phase transition at the point $\mu_\alpha=2\sec\alpha$. More precisely, whenever $\mu\leq \mu_\alpha$, the incoherent (unsynchronized) state $\rho_\infty\equiv\frac1{2\pi}$ is the unique solution. For $\mu>\mu_\alpha$, there exists an additional coherent traveling-wave solution 
\[
\rho(t,\theta)=\rho_{EMvM}(\theta-kt;\theta_0,k,r),
\]
which is unique up to the arbitrary choice of $\theta_0\in[0,2\pi)$, and where $k$ and $r$ are solutions of equation (\ref{OP3}). 

For $\alpha=0$, equation~\eqref{trwaves} undergoes a phase transition at the point $\mu=2$ to a coherent equilibrium steady-state $\rho(\theta)=\rho_{vM}(\theta;\theta_0,r)$, which is unique up to the arbitrary choice of $\theta_0\in[0,2\pi)$. 

For $\alpha=\pm\pi/2$, the incoherent state is unique regardless of $\mu$.    
\end{theorem}
\begin{remark}
    The above results apply to the case $|\alpha|\leq\pi/2$. To understand what happens when $\pi\geq|\alpha|>\pi/2$, note that the problem is invariant under the change of variables $(\mu,\alpha)\mapsto(-\mu,\pi-\alpha)$. Therefore, we obtain analogous results, with the difference that the bifurcation parameter is $-\mu$, and consequently the supercritical pitchfork bifurcation/phase transition in $-\mu$ becomes subcritical in $\mu$.  
\end{remark}
\section{Isogonal monotonicity of the mean resultant length and the resultant disk}
\subsection{Representation by power series with positive coefficients}
To prove the main results, we will need a series of lemmata.  
\begin{lemma}
For $p=0,1,2,\dots$, we define,  
    \be\label{rho}
    A _{p}(k)=
     \int_{-\pi}^{\pi}\cos^{2p}\left(\th/2\right)\exp(k\theta)\ d\theta.   
    \ee
Then we have a recursion relation, 
    \[
   A_{p}(k)=\frac{2p(2p-1)}{4(p^2+k^2)}A_{p-1}(k), 
    \]
and hence 
    \be\label{ap}
   A_{p}(k)=\frac{(2p)!}{4^p\Pi_{n=1}^p(n^2+k^2)}\frac{2\sinh(k\pi)}{k}.  
    \ee
\end{lemma}
\proof It follows by integration by parts.   \proofend
\begin{lemma}\label{parity}
    \begin{enumerate}
  \item We have, 
\begin{align}
\mathbb{I}_0(k,r) =C_0(k,r)&=\sum_{n=-\infty}^\infty\frac{k^2}{k^2+n^2}(-1)^nI_n(r)^2\label{C0}\\
&=1+\sum_{n=-\infty}^\infty\frac{n^2}{k^2+n^2}(-1)^{n+1}I_n(r)^2 \label{C00}
\end{align}
and, in particular, $\mathbb{I}_0(k,0)=1$ and $\mathbb{I}_0(0,r)=I_0(r)^2$. Furthermore, we have the following power series expansion:
\be\label{PSE}
C_0(k,r)=\sum_{p=0}^\infty a_p(k)\left(\frac{r}{2}\right)^{2p}, 
\ee
with $a_0(k)=1$ and 
\begin{align}
a_p(k)&=\frac1{p!^2}\sum_{n=-p}^p (-1)^{n+1}\frac{n^2}{n^2+k^2}{{2p}\choose{p+n}}\label{akp} \\
&=\frac{(2p)!}{p!^2\Pi_{n=1}^p(n^2+k^2)},\ p=1,2,\dots\label{akp1}
\end{align} 
\item We further have, 
\begin{align}
C_1(k,r)&=\sum_{n=-\infty}^{\infty}\frac{k^2}{k^2+n^2}(-1)^nI_n(r)I_n'(r)\label{C1}\\
&=\frac12\sum_{p=1}^\infty pa_p(k)\left(\frac{r}{2}\right)^{2p-1}, \label{C11}
\end{align} 
\begin{align}
S_1(k,r)&=-\frac kr \sum_{n=-\infty}^{\infty}\frac{n^2}{k^2+n^2}(-1)^nI_n(r)^2\label{S1}\\
&=\frac{1}2\sum_{p=1}^\infty ka_p(k)\left(\frac{r}{2}\right)^{2p-1},\label{S11}
\end{align}
and consequently, 
\[
\mathbb{I}_1(k,r)=\frac12\sum_{p=1}^\infty (p+ik)a_p(k)\left(\frac{r}{2}\right)^{2p-1}. 
\]
Furthermore, 
\be\label{T1}
T_1(k,r)=\frac{\sum_{p=1}^\infty ka_{p}(k)\left(\frac{r}{2}\right)^{2p}}{\sum_{p=1}^\infty p a_{p}(k)\left(\frac{r}{2}\right)^{2p}} 
\ee
 and 
\begin{align}\label{CC1}
\mathcal{C}(k,r)&=\frac{\sum_{p=0}^\infty (p+1)a_{p+1}(k)\left(\frac{r}{2}\right)^{2p}}{2\sum_{p=0}^\infty a_{p}(k)\left(\frac{r}{2}\right)^{2p}} \\
&=\frac{\sum_{p=0}^\infty (p+1)a_{p+1}(k)\left(\frac{r}{2}\right)^{2p}}{\sum_{p=0}^\infty\frac{p+1}{2p+1}((p+1)^2+k^2) a_{p+1}(k)\left(\frac{r}{2}\right)^{2p}}.\nonumber
\end{align}
\item We have:
\begin{enumerate}
    \item $\mathbb{I}_0(k,r)=C_0(k,r)\geq1$ with the equality holding for $r=0$; $C_0$ is even with respect to both $k$ and $r$.  
    \item $C_1(k,r)\geq0$ ($\leq0$) for $r\geq0$ ($r\leq0$), with equality holding for $r=0$; $C_1$ is even with respect to $k$ and odd with respect to $r$. 
    \item $S_1(k,r)\geq0$ ($\leq0$) for $kr\geq0$ ($kr\leq0$), with the equality holding for $kr=0$. $S_1$ is odd with respect to both $k$ and $r$.  
    \item $T_1$ is odd with respect to $k$ and even with respect to $r$. 
\end{enumerate}
    \end{enumerate}
\end{lemma}
\proof
\begin{enumerate}
    \item Since $\mathbb{I}_0(k,r)=\frac1{2\pi}\int_{0}^{2\pi}\chi_{k,r}(\phi)\ d\phi$, \eqref{chi} and~\eqref{psi} yield (\ref{C0}). Relation (\ref{one}) then yields (\ref{C00}). We then use (\ref{quadratic}) to obtain the power series expansion (\ref{PSE}), where the coefficients $a_p(k)$ are given by (\ref{akp}). To obtain formula~\eqref{akp1}, we first introduce 
    \[
    b_p(k)=\frac1{p!^2}\sum_{n=-p}^p (-1)^{n+1}\frac{k^2}{n^2+k^2}{{2p}\choose{p+n}}, 
    \]
    and then observe that by virtue of the Binomial Theorem, 
    \[
    a_p(k)+b_p(k)=\frac1{p!^2}\sum_{n=-p}^p (-1)^{n+1}{{2p}\choose{p+n}}=0.  
    \]
    Recall the integral representation of the binomial coefficients, 
    \[
    {{2p}\choose{p+n}}=\frac1{2\pi i}\oint_{|z|=1}\frac{(1+z)^{2p}}{z^{p+n+1}}\ dz.  
    \]
     Hence 
    \[
    b_p(k)=\frac1{p!^22\pi i}\oint_{|z|=1}\frac{F_{k,p}(z)}{z^{p+1}}, 
    \]
    where 
    \[
    F_{k,p}(z)=(1+z)^{2p}\sum_{n=-p}^p (-1)^{n+1}\frac{k^2}{n^2+k^2}z^{n}.
    \]
    By virtue of Laurent's theorem, $p!^2b_p(k)$ is the $p$\textsuperscript{th} coefficient for the rational function $F_{k,p}(z)$, which remain the same in $F_{k,p'}(z)$ for $p'>p$, as well as in the Laurent series, 
    \[
    F_{k,\infty}(z)=(1+z)^{2p}\sum_{n=-\infty}^\infty (-1)^{n+1}\frac{k^2}{n^2+k^2}z^{n}.
    \]
    Therefore, recalling the Fourier expansion 
    \[
    \sum_{n=-\infty}^{\infty}\frac{k^2}{n^2+k^2}e^{in\th}=\pi k\frac{\cosh(k(\pi-|\theta|))}{\sinh(k\pi)}, \th\in[-\pi,\pi], 
    \]
    we obtain 
    \beas
    b_p(k)&=&-\frac1{p!^22\pi}\int_{-\pi}^\pi(1+e^{i\th})^{2p}e^{-ip\th}\sum_{n\in\mathbb{Z}}\frac{k^2}{n^2+k^2}e^{n(\pi+\th)i}\ d\th\\
    &=&-\frac{2^{2p-1}}{p!^2}\int_{-\pi}^{\pi}\cos^{2p}\left(\th/2\right)\frac{k\cosh(k\theta)}{\sinh(k\pi)}\ d\theta\\
    &=&-\frac{2^{2p-1}}{p!^2}\frac{k}{\sinh(k\pi)}A_p(k). 
    \eeas
    Therefore, 
    \[
    a_p(k)=-b_p(k)=\frac{2^{2p-1}}{p!^2}\frac{k}{\sinh(k\pi)}A_p(k), 
    \]
    and the expression~\eqref{akp1} follows by~\eqref{ap}. 
    \item The relations (\ref{C1}) and (\ref{S1}) follow by taking the Fourier transform in (\ref{chi}), then applying trigonometric addition formulas and recursion formulas (\ref{rec1}) and (\ref{rec2}), respectively. In view of~\eqref{C0}, we obtain 
    \[
    C_1(k,r)=\frac12 \partial_r \mathbb{I}_0(k,r)
    \]
    and 
    \[
    S_1(k,r)=\frac kr\left(\mathbb{I}_0(k,r)-1\right),
    \]
    so~\eqref{C00} yields~\eqref{C11} and~\eqref{S11}. Relations~\eqref{T1} and~\eqref{CC1} follow then immediately in combination with the recursion formula for $a_p(k)$ for the latter.   
    \item follows from (1) and (2). 
\end{enumerate}  
\proofend
\begin{remark}
The two representations of the moment coefficients play complementary and non-overlapping
roles. The identity $\mathbb{I}_n(k,r)=|\Gamma(1+ik)|^2\,I_{ik}(r)\,I_{n-ik}(r)$ of
Proposition~\ref{prop:cont} is the structural one: it furnishes the three-term
recursion, the analytic continuation, and the identification of EMvM with the
classical Bessel family at imaginary order. It is, however, analytically opaque
to positivity, since $I_{n-ik}(r)$ carries an oscillatory phase (through
$(r/2)^{-ik}$) and no sign structure is visible in it. All monotonicity arguments
below instead operate on the \emph{real power-series} representation
\[
   \mathbb{I}_0(k,r)=\sum_{p=0}^\infty a_p(k)\Big(\frac r2\Big)^{2p},
   \qquad a_p(k)=\frac{(2p)!}{p!^2\prod_{n=1}^p(n^2+k^2)}>0,
\]
in which the frustration $k$ enters only through the strictly positive factor
$\prod_{n=1}^p(n^2+k^2)$ in the denominator. It is this positivity of the
coefficients that is preserved under
the exponential skewing. The complex-order form is the bridge to the classical theory; positivity is a
power-series phenomenon, and the proofs of the subsequent results exploit the latter.
\end{remark}

Monotonicity arguments largly rest on the following lemma due to Biernacki and Krzyż (see~\cite{BK, YangChuWang}). 
\begin{lemma}\label{krzyz}
    Let $\alpha_p>0$, and $\beta_p>0$ for $p=0,1,2,\dots$, and let the power series $f(r)=\sum_{p=0}^\infty \alpha_{p}x^{p}$ and $g(r)=\sum_{p=0}^\infty\beta_{p}x^{p}$ converge for $x>0$. If the quotient sequence 
    \[
    \delta_p=\frac{\alpha_p}{\beta_p}, \ p=0,1,2,\dots
    \]
    is strictly increasing (decreasing), then the function $W(x)=\frac{f(x)}{g(x)}$
    is also strictly increasing (decreasing) on $\mathbb{R}_{>0}$. 
\end{lemma}  
\begin{lemma}\label{first}
\begin{enumerate}
    \item The function 
    \[
    \mathbb{R}_{>0}\owns r\mapsto T_1(k,r)\in(0,k)
    \]
    is a strictly decreasing bijection for $k>0$. Moreover, for $k>0$,  
    \be\label{lim1}
    \lim_{r\to\infty} \frac{r}{k} T_1(k,r)=1.
    \ee
    \item The function 
    \[
    \mathbb{R}_{>0}\owns k\mapsto \frac1k T_1(k,r)\in\left(\frac{I_0(r)^2-1}{rI_0(r)I_1(r)},1\right)
    \]
    is a strictly increasing bijection for $r>0$. 
    In particular, 
    \be\label{lim2}
    \lim_{k\to\infty}\frac1k T_1(k,r)=1.
    \ee
    \item 
    For $k\in\mathbb{R}$ and $r>0$, $\partial_{k} T_1>0$. The function $\mathbb{R}\owns k\mapsto T_1(k,r)\in\mathbb{R}$ is a strictly increasing bijection for $r>0$. 
    \item The function $\mathcal{C}(k,r)$ is strictly decreasing in $k>0$ for a fixed $r>0$, but it is not necessarily monotonic in $r>0$ for a fixed $k\not=0$. 
    \end{enumerate} 
    \end{lemma}
\proof 
\begin{enumerate}
    \item The fact that $ T_1(k,r)$ is strictly decreasing in $r>0$ for $k>0$ follows by virtue of Lemma~\ref{krzyz} applied to the quotient of power series (\ref{T1}). Further, it can be easily verified that 
    \begin{align*}
        \lim_{r\to0+}\mathcal{C}(k,r)&=\frac1{k^2+1}\\
        \lim_{r\to0+}\mathcal{S}(k,r)&=\frac k{k^2+1}\\
        \lim_{r\to0+}\mathcal{R}(k,r)&=\frac1{\sqrt{k^2+1}}\\
        \lim_{r\to0+} T_1(k,r)&=k. 
    \end{align*}
    The statement about the limit follows using the fact that $\lim_{r\to\infty}\frac{I_n(r)}{I_0(r)}=1$ and $\lim_{r\to\infty}\frac{I_n'(r)}{I_0(r)}=1$:  
 \begin{align*}
 \lim_{r\to\infty} \frac{r}{k} T_1(k,r)&=\lim_{r\to\infty} \frac{C_0(k,r)-1}{C_1(k,r)}\\
 &=\lim_{r\to\infty}\frac{\sum_{n\in\mathbb{Z}}\frac{(-1)^nk^2}{k^2+n^2}\frac{I_n(r)^2}{I_0(r)^2}-\frac1{I_0(r)^2}}{\sum_{n\in\mathbb{Z}}\frac{(-1)^nk^2}{k^2+n^2}\frac{I_n'(r)I_n(r)}{I_0(r)^2}}\\
 &=1.
\end{align*}

 \item We need to prove that for $r>0$,
 \[
 k\mapsto \frac1k T_1(k,r)=\frac{\sum_{p=1}^\infty a_p(k)\left(\frac{r}{2}\right)^{2p-1}}{\sum_{p=1}^\infty pa_p(k)\left(\frac{r}{2}\right)^{2p-1}}
 \]
 is strictly increasing. This follows as the derivative with respect to $k$ becomes 
 \[
      \partial_k\left(\frac1k T_1\right)(k,r)=\left(\sum_{p=1}^\infty c_{p}(k)\left(\frac{r}{2}\right)^{2p-1}\right)/\left(\sum_{p=1}^\infty pa_p(k)\left(\frac{r}{2}\right)^{2p-1}\right)^2,
 \]
 where 
 \begin{align*}
     c_{p}(k)&=\sum_{n=1}^{p-1}(p-n) (\partial_k a_{n}a_{p-n}- a_{n}\partial_k a_{p-n})(k)\\
        &=\sum_{n=1}^{[\frac{p-1}{2}]}(p-2n) a_{n}(k)a_{p-n}(k)\left(\partial_k(\log(a_{n}))- \partial_k(\log(a_{p-n}))\right)(k)\\
        &=\sum_{n=1}^{[\frac{p-1}{2}]}(p-2n) a_{n}(k)a_{p-n}(k)\sum_{j=n+1}^{p-n}\frac{2k}{j^2+k^2}>0,\ k>0. 
 \end{align*}
 The statement about the limit is obvious.
 \item By the computation above, 
 \[
 \partial_k T_1=\frac1k T_1+k \partial_k\left(\frac1k T_1\right)>0,\ k>0.  
 \]
 The bijection claim follows from~\eqref{lim2}. 
 \item The first statement follows directly from (2). It can be easily verified that for a fixed $k>0$, $\partial_r \mathcal{C}$ does not have a sign. This is also consistent with the fact that Lemma~\ref{krzyz} does not apply as    
 \[
 \frac{k^2+(p+1)^2}{2p+1}
 \]
 is not necessarily monotonic in $p$ except for $k=0$. 
 \end{enumerate}
 \proofend
\subsection{Isogonal coordinates}\label{ssec:isogone}
The following lemma introduces the change to isogonal coordinates. 
\begin{lemma}\label{isogone1}
The mapping 
\[
\mathbb{R}\times\mathbb{R}_+\owns (k,r)\mapsto(\kappa,r)=( T_1(k,r),r)\in\mathbb{R}\times\mathbb{R}_+
\]
is a diffeomorphism. Let 
\[
\mathbb{R}\times\mathbb{R}_+\owns (\kappa,r)\mapsto(k,r)=(K(\kappa,r),r)\in\mathbb{R}\times\mathbb{R}_+
\]
be its inverse. In particular, 
\[
\mathbb{R}\times\mathbb{R}_+=\cup_{\kappa\in\mathbb{R}}\Gamma_\kappa,  
\]
where 
\[
\Gamma_\kappa=\left\{(k,r)\in\mathbb{R}_+^2\ \Big|\ T_1(k,r)=\kappa\right\}=\left\{(K(\kappa, r), r)\ \Big|\ r>0\right\}. 
\]
Furthermore, $K(\kappa,r)$ has the sign of $\kappa$, and it is odd with respect to $\kappa$ with with $K(0,\cdot)\equiv0$; for $\kappa>0$,
$K(\kappa,\cdot):(0,\infty)\to(\kappa,\infty)$ is a strictly increasing bijection, while for $\kappa<0$,
$K(\kappa,\cdot):(0,\infty)\to(-\infty,\kappa)$ is a strictly decreasing bijection. 
\end{lemma} 
\proof 
For a fixed $r>0$ the map $k\mapsto T_1(k,r)$ is a strictly increasing bijection of $\mathbb{R}$ onto
$\mathbb{R}$ with $\partial_kT_1>0$ (Lemma~\ref{first}(3)); by the implicit function theorem the inverse
$(\kappa,r)\mapsto(K(\kappa,r),r)$ exists and is $C^1$, so the mapping is a diffeomorphism and
$\mathbb{R}\times\mathbb{R}_+$ is partitioned into the level sets $\Gamma_\kappa$. The ranges of
$K(\kappa,\cdot)$ are read off from Lemma~\ref{first}(1)--(2). Since $T_1(\cdot,r)$ is odd in $k$ (Lemma~\ref{parity}), $T_1(K,r)=\kappa$ forces $T_1(-K,r)=-\kappa$, that is,
\begin{equation}\label{eq:reflect}
   K(-\kappa,r)=-K(\kappa,r),\qquad \kappa\in\mathbb{R},\ r>0 .
\end{equation}

Throughout the rest of the proof we fix $\kappa\ge0$, and write for simplicity $K=K(\kappa,r)$,
$x=(r/2)^2$. Let $\rho$ be the probability measure on $\mathbb{Z}_{\ge0}$ with weights
\begin{equation}\label{eq:rhoweights}
   \rho_p\propto w^\rho_p:=a_{p+1}(K)\,x^{p},
\end{equation}
which has $\rho_p>0$ for all $p$ and, since $a_{p+1}(K)=O(4^p/(p!)^2)$ from~\eqref{akp1},
finite moments of every order. Expectations $\E$, variances $\Var$, and covariances $\Cov$ without a
subscript are taken against $\rho$. We further introduce the observable $A_p=p+1$, and the size--biased law $\mu$ with
$\mu_p\propto w^{\mu}_p:=A_p\,a_{p+1}(K)\,x^p$, so that $\mu_p=A_p\rho_p/\E[A]$. It will be convenient to
work with the \emph{unnormalized} weights of the two laws, $w^{\nu}_p$ for $\nu\in\{\rho,\mu\}$, 
so that $\nu_p=w^\nu_p/\sum_q w^\nu_q$. Equation~\eqref{T1} becomes 
\[
   T_1=\frac{S_1}{C_1}
   =\frac{K\sum_{p=0}^\infty a_{p+1}(K)x^p}{\sum_{p=0}^\infty(p+1)a_{p+1}(K)x^p}=\frac{K}{\E[A]},
\]
and the isogone constraint $T_1(K,r)=\kappa$ now reads
\begin{equation}\label{eq:isoconstraint}
   K=\kappa\,\E[A]:
\end{equation}
the bare frustration equals the isogone invariant times the size--bias normalizer (that also depends on $K$, and as such is not to be misunderstood for an explicit formula for $K$). 

The monotonicity of $K$ in $r$ will follow directly from the following lemma.  
\proofend
\begin{lemma}[isogonal score and slope]\label{score}
Fix an isogone $\Gamma_\kappa$, $\kappa>0$. We difine, 
\begin{equation}\label{eq:sigmadef}
   s_p:=-\partial_K\log a_{p+1}(K)=\sum_{n=1}^{p+1}\frac{2K}{n^2+K^2},\qquad
   \sigma_p:=K\,s_p=\sum_{n=1}^{p+1}\frac{2K^2}{n^2+K^2}. 
\end{equation}
\begin{enumerate}
\item[(i)] Define the ``score", 
\[
   \psi_p:=2p-\Big(r\frac{dK}{dr}\Big)\,s_p. 
\]
Then for \emph{both} $\nu\in\{\rho,\mu\}$,
\[
   r\frac{d}{dr}\log w^\nu_p=\psi_p,
\]
and consequently, for any observable $h$ of finite second moment,
\[
   r\frac{d}{dr}\E_\nu[h]=\Cov_\nu(h,\psi).
\]
\item[(ii)] Let 
\[
\Lambda=\Lambda(\kappa,r):=\dfrac{d\log K}{d\log r}\Big|_{\Gamma_\kappa}=\frac{r\partial_rK(\kappa,r)}{K(\kappa,r)}
\]
so that $r\,d K/dr=K\Lambda$. Then the score takes the form $\psi_p=2p-\Lambda\,\sigma_p$, and defining 
\begin{equation}\label{eq:Wpos}
   W:=\Cov(A,\sigma)=\Cov(p,\sigma) 
\end{equation}
we have 
\[
   \Lambda=\frac{\Cov(A,2p)}{\E[A]+W}
   =\frac{2\Var(p)}{\E[p]+1+W}>0.
\]
\end{enumerate}
If $\kappa=0$, then $K\equiv0$, $s\equiv\sigma\equiv0$ and $\psi_p=2p$. If $\kappa>0$, then $K>0$ on $\Gamma_\kappa$ and, 
$r\,dK/dr=K\Lambda>0$; hence $r\mapsto K(\kappa,r)$ is strictly increasing. For
$\kappa<0$ the reflection~\eqref{eq:reflect} turns this into strict decrease.
\end{lemma}
\proof
The two weights satisfy $\log w^{\mu}_p=\log (p+1)+\log w^\rho_p$, so it suffices to compute the $r$--derivative for $\log w^\rho_p=\log a_{p+1}(K)+ p\log x$. Since $r\frac{d}{dr}\log x^p=2p$, and $\partial_K\log a_{p+1}(K)=-s(p)$ by~\eqref{eq:sigmadef},
the chain rule along $\Gamma_\kappa$ gives
\[
   r\frac{d}{dr}\log a_{p+1}(K)=\big(\partial_K\log a_{p+1}\big)\cdot r\frac{dK}{dr}
   =-r\frac{dK}{dr}\,s_p.
\]
(Since $\rho_p=O(4^p/(p!)^2\cdot x^p)$ decays super--exponentially, the series
for $\E_\nu[h]$ and its termwise $r$--derivative converge locally uniformly, so differentiation under the
sum is justified.) Writing $\E_\nu[h]=\sum_p h_pw^\nu_p/\sum_p w^\nu_p$ and differentiating, we obrtain 
$r\frac{d}{dr}\E_\nu[h]=\E_\nu[h\psi]-\E_\nu[h]\,\E_\nu[\psi]=\Cov_\nu(h,\psi)$. This proves (i).

For (ii), substituting $r\,dK/dr=K\Lambda$ and $K\,s(p)=\sigma(p)$ into the
score of (i) yields $\psi_p=2p-\Lambda\sigma(p)$. Applying $r\frac{d}{dr}\log(\cdot)$ to
\eqref{eq:isoconstraint}, with $\kappa$ constant and $r\frac{d}{dr}\log\E[A]=\Cov(A,\psi)/\E[A]$ by (i),
\[
   \Lambda=\frac{\Cov(A,\psi)}{\E[A]}=\frac{\Cov(A,2p)-\Lambda\Cov(A,\sigma)}{\E[A]}.
\]
Solving this linear relation for $\Lambda$ gives the first form; $A=p+1$ gives the second. For the sign,
note that $\sigma$ is nondecreasing, since $\Delta\sigma(p)=2K^2/\big((p+2)^2+K^2\big)\ge0$; as $p$ is
increasing as well, Chebyshev's sum inequality gives $W\geq0$. 
Hence the denominator satisfies $\E[p]+1+W\ge1>0$, while $\Var(p)>0$ because
$\rho_p>0$ for all $p$; therefore $\Lambda>0$. 
\proofend
While abusing notation slightly, we consider now the mean resultant length in isogonal coordinates $(\kappa,r)$,  
\be\label{Ctilde}
\mathcal{R}(\kappa,r)=\mathcal{R}(K(\kappa,r),r),\ \kappa\in\mathbb{R}, r>0,
\ee
(similarly $\mathcal{C}(\kappa,r)$, $\mathcal{S}(\kappa,r)$). Recall that $T_1(\kappa,r)=\frac{\mathcal{S}(\kappa,r)}{\mathcal{C}(\kappa,r)}=\kappa$ is constant, so $\mathcal{S}(\kappa,r)=\kappa\mathcal{C}(\kappa,r)$, and ${\mathcal R}(\kappa,r)=\sqrt{1+\kappa^2}\,\mathcal{C}(\kappa,r)$.  
Substituting
\eqref{eq:isoconstraint} into the identity~\eqref{CC1} for $\mathcal{C}$ yields 
\begin{equation}\label{eq:invCsplit}
   \frac1{\mathcal C(\kappa,r)}=F_1+\kappa^2F_2,\qquad
   F_1=\E_{\mu}[b_1],\qquad F_2=\E[A]\,\E[b_0],
\end{equation}
where the observables are the size--bias weights $A_p=p+1$, $B_p=\frac{p+1}{2p+1}$, and the one--parameter family
\begin{equation}\label{eq:obsfamily}
   b_j(p):=A_p^{\,j}\,B_p=\frac{(p+1)^{j+1}}{2p+1},\qquad j=0,1,2. 
\end{equation}
In this notation the first term of~\eqref{eq:invCsplit} is the size--biased mean $\E_\mu[b_1]$ and the
second the product of two $\rho$--means $\E[A]\,\E[b_0]$. 
\begin{lemma}\label{isogone2}
\begin{enumerate}
\item[(i)] For $\kappa\in\mathbb{R}$, the map
\[
(0,\infty)\owns r\mapsto \mathcal{R}(\kappa,r)\in(0,1/\sqrt{1+\kappa^2})
\]
is a strictly decreasing bijection.
\item[(ii)] For every $k\in\mathbb{R}$ and $r>0$ one has $\Re(1/\mathcal{M})>1$; equivalently
$|\mathcal{M}|^2<\Re\mathcal{M}$, equivalently $\mathcal{R}^2<\mathcal{C}$, equivalently
$|\mathcal{M}-\tfrac12|<\tfrac12$, equivalently $\mathcal{R}<\cos\mathcal{A}$. That is,
$\mathcal{M}(k,r)$ lies in the open disk $\mathcal{W}=\{z\in\mathbb{C}:|z-\tfrac12|<\tfrac12\}$.
\end{enumerate}
\end{lemma} 
\begin{remark}
Since the weights $a_{p+1}(K)$ see $K$ only through $K^2$, the measures $\rho$ and $\mu$, and every observable built from them, are therefore unchanged under $\kappa\mapsto-\kappa$. By
\eqref{eq:invCsplit} the function ${\mathcal C}$, and therefore
${\mathcal R}=\sqrt{1+\kappa^2}\,\mathcal{C}$, depend on $\kappa$ only through
$\kappa^2$. It is thus enough to prove Lemma~\ref{isogone2} for $\kappa\ge0$, the case $\kappa<0$
following by this reflection; \emph{we assume $\kappa\ge0$ from now on}, so that $K\ge0$ on
$\Gamma_\kappa$, with $K>0$ precisely when $\kappa>0$ (by~\eqref{eq:isoconstraint} below).
\end{remark}
\subsection{Two covariance inequalities}
For the following lemma, let $\rho$ be a probability
measure on $\mathbb{Z}_{\ge0}$, and $\E$, $\Var$, $\Cov$ denote expectation, variance and covariance
against $\rho$; $p$ denotes the identity observable $p\mapsto p$, so that $\E[p]$, $\Var(p)$ and
$\Cov(h,p)$ have their usual meaning. An \emph{observable} $h=\{h_p\}_{p\ge0}$ is
\emph{convex} (resp.\ \emph{concave}) if its second difference $\Delta^2h(j):=h_{j+2}-2h_{j+1}+h_j$ is
$\ge0$ (resp.\ $\le0$) for all $j\ge0$. The lemma concerns the sign for a covariance determinant. For observables $f,h$ of finite second moment set
\begin{equation}\label{eq:Ddef}
   D(f,h):=\Var(p)\,\Cov(f,h)-\Cov(p,h)\,\Cov(f,p)
   =\det\!\begin{pmatrix}\Cov(p,p)&\Cov(p,h)\\[2pt]\Cov(f,p)&\Cov(f,h)\end{pmatrix}.
\end{equation}
\begin{lemma}\label{det}
Let $\rho$ be a probability measure on $\mathbb{Z}_{\ge0}$ with finite second moment and be non-degenerate
$($i.e.\ $\Var(p)>0)$, and let $f,h$ be observables of finite second moment against $\rho$. If $f$ and $h$ are convex, then $D(f,h)\ge0$. Moreover, when $\rho_p>0$ for all $p$ the inequality is \emph{strict} whenever $\Delta^2f(0)>0$ and $\Delta^2h(0)>0$.
\end{lemma}
\proof
This is the one–dimensional determinantal covariance inequality of Bonnefont, Hillion and Saumard
\cite[Theorem~1.2(1), eq.~(3.11)]{BHS}, applied to the identity observable $p$ and to the two convex
observables $f$ and $h$. Their result holds for a general probability measure with finite second
moment, of which our discrete $\rho$ is a special case.
\proofend

\begin{remark}\label{rem:strictdet}
The mechanism behind Lemma~\ref{det} is classical: via Hoeffding's covariance identity the quantity
$D(f,h)$ is itself a covariance against the bivariate measure with the
totally positive kernel $\Cov(\mathbf 1_{\{p>v\}},\mathbf 1_{\{p>w\}})$, whose sign follows from the FKG
inequality~\cite{FKG} applied to the increasing observables $\Delta f$, $\Delta h$ (see
\cite[\S3]{BHS} and, for the total positivity, Karlin~\cite{Karlin}; the Gaussian antecedents are
Hu~\cite{Hu} and Harg\'e~\cite{Harge}). 
\end{remark}

\begin{lemma}\label{krzyz1}
Let $\rho$ be a probability measure on $\mathbb{Z}_{\ge0}$ with finite third moment and be non-degenerate.
With $A_p=p+1$ and $B_p=\frac{p+1}{2p+1}$ as in~\eqref{eq:obsfamily},
\begin{equation}\label{eq:dom}
   \Cov(A,p)\,\E[B]+\Cov(B,p)\,\E[A]>0 .
\end{equation}
\end{lemma}
\proof
As $A_p$ increases and $B_p$ decreases, $\Cov(A,p)\ge0$ and $\Cov(B,p)\le0$. Writing the two
covariances as double sums against $\rho$ and symmetrizing, the left side of~\eqref{eq:dom} equals
$\sum_{p,q,s}\rho_p\rho_q\rho_s\,T(p,q,s)$ with
\[
   T(p,q,s)=\frac1{|S_3|}\!\!\sum_{\sigma\in S_3}\!\!
   \big[(A_{\sigma p}-A_{\sigma q})(\sigma p-\sigma q)B_{\sigma s}
   +(B_{\sigma p}-B_{\sigma q})(\sigma p-\sigma q)A_{\sigma s}\big].
\]
Here $S_3$ is the group of the six permutations of the three indices $(p,q,s)$, and
$\sigma p,\sigma q,\sigma s$ denote the permuted indices ($|S_3|=6$). The symmetrization is
justified by the fact that the weight product $\rho_p\rho_q\rho_s$ is symmetric in $(p,q,s)$, so summing
the raw (asymmetric) summand over all $p,q,s$ equals summing its $S_3$-average $T$, which is a symmetric
function of $(p,q,s)$. Since the weights are nonnegative, and on the diagonal $p=q=s$ one has $T=0$
(each summand carries a factor $\sigma p-\sigma q$), it suffices to show $T>0$ off the diagonal.

Using $A_a-A_b=a-b$ and $B_a-B_b=-\frac{a-b}{(2a+1)(2b+1)}$, each summand collapses and combining
the six permutations gives the exact identity
\[
   T(p,q,s)=\frac{2\,N(p,q,s)}{3\,(2p+1)(2q+1)(2s+1)},
\]
with
\[
   N=(p{-}q)^2(s{+}1)(2pq{+}p{+}q{-}s)
   +(p{-}s)^2(q{+}1)(2ps{+}p{+}s{-}q)
   +(q{-}s)^2(p{+}1)(2qs{+}q{+}s{-}p).
\]
Because of the symmetry, we may assume $p\ge q\ge s\ge0$ and set
$s\ge0$, $d=q-s\ge0$, $e=p-q\ge0$. Then
\[
\begin{aligned}
   N=\;&4d^4s+8d^3es+8d^3s^2+6d^3s+8d^2e^2s+3d^2e^2+12d^2es^2+9d^2es\\
   &+4d^2s^3+6d^2s^2+2d^2s+4de^3s+3de^3+12de^2s^2+15de^2s+4de^2\\
   &+4des^3+6des^2+2des+4e^3s^2+6e^3s+2e^3+4e^2s^3+6e^2s^2+2e^2s,
\end{aligned}
\]
all coefficients nonnegative, so $N\ge0$; every monomial carries a factor $d$ or $e$, so $N=0$ only
when $p=q=s$. Since $\Var(p)>0$, the measure $\rho$ assigns positive mass to at least two distinct
points, hence to some off-diagonal triple $(p,q,s)$; there $T>0$, so the triple sum in~\eqref{eq:dom}
is strictly positive.
\proofend
\begin{proposition}\label{derivs}
Write $W$ for the quantity in~\eqref{eq:Wpos}. 
Let 
\[
   D_0:=\E[b_0]\Cov(A,2p)+\E[A]\Cov(b_0,2p),\qquad
   D_1:=\E[b_0]\Cov(A,\sigma)+\E[A]\Cov(b_0,\sigma),
\]
and 
\[
   D_0^\sharp:=\E[A]\Cov(b_2,2p)-\E[b_2]\Cov(A,2p),\qquad
   D_1^\sharp:=\E[A]\Cov(b_2,\sigma)-\E[b_2]\Cov(A,\sigma). 
\]
Then $D_0>0$ and $D_0^\sharp>0$. Moreover:
\begin{enumerate}
\item[(i)] Along any isogone $\Gamma_\kappa$,
\[
   r\frac{dF_1}{dr}=\Cov_{\mu}(b_1,\psi),\qquad
   r\frac{dF_2}{dr}=\E[b_0]\,\Cov(A,\psi)+\E[A]\,\Cov(b_0,\psi).
\]
\item[(ii)] On $\Gamma_0$ one has $\psi=2p$, so that
\[
   r\frac{dF_1}{dr}=\E[A]^{-2}D_0^\sharp>0\qquad\text{and}\qquad r\frac{dF_2}{dr}=D_0>0 .
\]
\item[(iii)] If $\kappa>0$, then with $D$ as in
\eqref{eq:Ddef},
\[
    r\frac{dF_1}{dr}=\frac1{\E[A](\E[A]+W)}\big(D_0^\sharp-2D(b_2,\sigma)\big), \qquad
   r\frac{dF_2}{dr}=\frac{\E[A]}{\E[A]+W}\big(D_0-2D(b_0,\sigma)\big).
\]
\end{enumerate}
\end{proposition}
\proof
The positivity $D_0>0$ is Lemma~\ref{krzyz1}, and $D_0^\sharp=\E[A]^2\Cov_{\mu}(b_1,2p)>0$ by the strict
Chebyshev sum inequality, since $b_1$ and $2p$ are strictly increasing and $\rho$ is non-degenerate.

(i) The first identity is Lemma~\ref{score}(i) with $h=b_1$, $\nu=\mu$. For the second,
$F_2=\E[A]\E[b_0]$, and the product rule together with Lemma~\ref{score}(i) (applied to $h=A$ and $h=b_0$,
$\nu=\rho$) give $r\,dF_2/dr=\E[b_0]\Cov(A,\psi)+\E[A]\Cov(b_0,\psi)$. 

(ii) On $\Gamma_0$ we have $\psi=2p$; substituting and using bilinearity of the covariance gives
$r\,dF_1/dr=\Cov_{\mu}(b_1,2p)=\E[A]^{-2}D_0^\sharp$ and $r\,dF_2/dr=D_0$, both positive by the
above.

(iii) Let $\kappa>0$. Substituting $\psi=2p-\Lambda\sigma$ (Lemma~\ref{score}(ii)) into the identities of
(i) and using bilinearity gives $r\,dF_2/dr=D_0-\Lambda D_1$ and
$r\,dF_1/dr=\E[A]^{-2}(D_0^\sharp-\Lambda D_1^\sharp)$. In each row the same cancellation occurs: expanding
the two $D$'s and using $A=p+1$ (so $\Cov(A,\cdot)=\Cov(p,\cdot)$ and $\Cov(\cdot,2p)=2\Cov(\cdot,p)$), the
leading terms cancel and the determinant~\eqref{eq:Ddef} emerges,
\beas
   \Cov(A,\sigma)\,D_0-\Cov(A,2p)\,D_1=-2\E[A]\,D(b_0,\sigma),\\
   \Cov(A,\sigma)\,D_0^\sharp-\Cov(A,2p)\,D_1^\sharp=-2\E[A]\,D(b_2,\sigma).
\eeas 
Since $\Lambda=\Cov(A,2p)/(\E[A]+W)$ and $W=\Cov(A,\sigma)$ (Lemma~\ref{score}(ii)),
\beas
   r\frac{dF_2}{dr}&=&D_0-\Lambda D_1\\
   &=&\frac{(\E[A]+W)D_0-\Cov(A,2p)D_1}{\E[A]+W}\\
   &=&\frac{\E[A]\,D_0+\big(\Cov(A,\sigma)D_0-\Cov(A,2p)D_1\big)}{\E[A]+W}\\
   &=&\frac{\E[A]}{\E[A]+W}\big(D_0-2D(b_0,\sigma)\big). 
\eeas
The identical manipulation with the second identity
yields the formula for $r\,dF_1/dr$.
\proofend
\proof[Proof of Lemma~\ref{isogone2}]
\emph{Part (i).} By~\eqref{eq:invCsplit}, $1/\mathcal C(\kappa,r)=F_1+\kappa^2F_2$, and
$\mathcal {R}(\kappa,r)=\sqrt{1+\kappa^2}\,\mathcal {C}$. It therefore suffices to show
$r\,dF_1/dr>0$ and $r\,dF_2/dr>0$ for every $r>0$.

Fix $\kappa>0$. The observables $b_0$ and $b_2$ are convex, $\sigma$ is
concave, and
\[
   \Delta^2b_0(0)=\tfrac4{15}>0,\qquad \Delta^2b_2(0)=\tfrac{16}{15}>0,\qquad
   \Delta^2\sigma(0)=\frac{2K^2}{9+K^2}-\frac{2K^2}{4+K^2}<0
\]
(the last strict because $K>0$ on $\Gamma_\kappa$ for $\kappa>0$). Indeed
$b_0=B_p=\tfrac12+\tfrac1{2(2p+1)}$ is a constant plus a convex term;
$b_2=\tfrac18\big((2p+1)^2+3(2p+1)+3\big)+\tfrac1{8(2p+1)}$ is a convex quadratic plus a convex term with
$\Delta^2b_2\to1$ bounded; and $\sigma$ has strictly decreasing increments $\tfrac{2K^2}{(p+2)^2+K^2}$,
hence is concave. Since $\rho_p>0$ for all $p$, Lemma~\ref{det} together with the strictness criterion
(the displayed $\Delta^2$ values give $\Delta^2b_0(0),\Delta^2b_2(0)>0$ and
$\Delta^2\sigma(0)<0$) gives
\[
   D(b_0,\sigma)<0\qquad\text{and}\qquad D(b_2,\sigma)<0.
\]
With $D_0>0$, $D_0^\sharp>0$ and the positive prefactors of Proposition~\ref{derivs}(ii),
\[
   r\frac{dF_2}{dr}=\frac{\E[A]}{\E[A]+W}\big(D_0-2D(b_0,\sigma)\big)>0,\qquad
   r\frac{dF_1}{dr}=\frac{D_0^\sharp-2D(b_2,\sigma)}{\E[A](\E[A]+W)}>0.
\]
The isogone $\Gamma_0$ is covered by Proposition~\ref{derivs}(i), which gives
$r\,dF_1/dr=\E[A]^{-2}D_0^\sharp>0$ and $r\,dF_2/dr=D_0>0$ there. Thus for both $\kappa>0$ and
$\kappa=0$, and hence by the reflection~\eqref{eq:reflect} for every $\kappa\in\mathbb{R}$, $\mathcal{C}$ and hence $\mathcal {R}=\sqrt{1+\kappa^2}\,\mathcal{C}$ are strictly
decreasing in $r$ along $\Gamma_\kappa$. The stated range and the bijection claim follow from the
boundary values $\lim_{r\to0+}\mathcal{R}(\kappa,r)=1/\sqrt{1+\kappa^2}$ (Lemma~\ref{first}(1)) and
$\lim_{r\to\infty}\mathcal{R}(\kappa,r)=0$.

\emph{Part (ii).} Fix $(k,r)$ with $r>0$; by the
decomposition~\eqref{eq:invCsplit}, 
\[
   \Re\!\Big(\frac1{\mathcal{M}}\Big)=\frac{1/\mathcal{C}}{1+\kappa^2}=\frac{F_1+\kappa^2F_2}{1+\kappa^2}.
\]
Thus $\Re(1/\mathcal{M})>1$ is equivalent to $(F_1-1)+\kappa^2(F_2-1)>0$, and it suffices to show
$F_1\ge1$ and $F_2\ge1$, strictly for $r>0$. The observable satisfies
$b_1(p)=(p+1)^2/(2p+1)\ge1$, since $(p+1)^2-(2p+1)=p^2\ge0$ (equality only at $p=0$); hence
$F_1=\E_{\mu}[b_1]\ge1$, strictly for $r>0$ because then $\mu$ is non-degenerate. Since
$A_p=p+1$ is increasing and $B_p=(p+1)/(2p+1)$ decreasing, Chebyshev's sum inequality gives
$\Cov(A,b_0)\le0$, i.e.\ $\E[b_1]=\E[A\,b_0]\le\E[A]\E[b_0]=F_2$; with $\E[b_1]\ge1$ this gives
$F_2\ge1$, strictly for $r>0$. Hence $\Re(1/\mathcal{M})>1$, i.e.\ $\mathcal{M}$ lies in the open disk
$\{|z-\tfrac12|<\tfrac12\}$. The stated equivalences are the identities
$\Re(1/\mathcal{M})>1\Leftrightarrow|\mathcal{M}|^2<\Re\mathcal{M}\Leftrightarrow\mathcal{R}^2<\mathcal{C}
\Leftrightarrow|\mathcal{M}-\tfrac12|<\tfrac12\Leftrightarrow\mathcal{R}<\cos\mathcal{A}$.
\proofend

The proof of Theorem~\ref{main0} is an immediate consequence of the following 
\begin{lemma}\label{wedge}
Let 
\[
\mathcal{V}=\left\{z=re^{i\beta}\ \Big|\ \beta\in(-\pi/2,\pi/2)\right\}
\]
denote the open half-plane $\mathbb{R}_{>0}\times\mathbb{R}$, and let 
\beas
\mathcal{W}&=&\left\{z=Re^{i\alpha}\ \Big|\ \alpha\in(-\pi/2,\pi/2),\ 0<R<\cos\alpha\right\}\\
 &=&\left\{z:|z-\tfrac12|<\tfrac12\right\}\\
 &=&\left\{z:\Re(1/z)>1\right\}
\eeas
denote the disk centered at $1/2$ of radius $1/2$. 
Then the normalized mean resultant map
\[
 \mathcal{V}\owns re^{i\beta}\mapsto\mathcal{M}(\tan\beta,r)=Re^{i\alpha}\in \mathcal{W}
\]
is strictly decreasing along the isogones and globally invertible, with radially strictly
decreasing inverse. It maps $\mathcal{V}$ diffeomorphically onto $\mathcal{W}$, mapping each isogone monotonically onto an open 
chord of disk $\mathcal{W}$ between points $z=0$ and $z=\cos\alpha e^{i\alpha}$.
\end{lemma}
\proof
The three descriptions of $\mathcal{W}$ agree: for $z=Re^{i\alpha}=x+iy$ with $R>0$, the bound $R<\cos\alpha$
reads $R^2<R\cos\alpha$, i.e.\ $x^2+y^2<x$, i.e.\ $|z-\tfrac12|<\tfrac12$; and
$|z-\tfrac12|<\tfrac12$ is $|z|^2<\Re z$, equivalently $\Re(1/z)=\Re z/|z|^2>1$.

By Lemma~\ref{isogone1}, setting $\kappa=\tan\alpha$, $\alpha\in(-\pi/2,\pi/2)$, the isogones of $\mathcal{M}$ are given by 
\[
\Gamma_{\tan\alpha}=\left\{(k,r)\in\mathbb{R}_+^2\ \Big|\ T_1(k,r)=\tan\alpha\right\}=\left\{(K(\tan\alpha, r), r)\ \Big|\ r>0\right\},
\]
and by Lemma~\ref{isogone2}(i) the normalized mean resultant length along these isogones,
\[
(0,\infty)\owns r\mapsto \mathcal{R}(\tan\alpha,r)\in(0,\cos\alpha)
\]
is a strictly decreasing bijection. Hence each isogone maps bijectively onto the open chord 
$\{Re^{i\alpha}:0<R<\cos\alpha\}$; as $\alpha$ ranges over $(-\pi/2,\pi/2)$ these chords are disjoint with
union $\mathcal{W}$, so $\mathcal{M}$ is a bijection of $\mathcal{V}$ onto $\mathcal{W}$ (injective by the global invertibility of
Lemma~\ref{isogone1}). It is moreover a diffeomorphism. Indeed, $\mathcal{M}$ factors as
\[
   re^{i\beta}\ \mapsto\ (k,r)=(\tan\beta,r)\ \mapsto\ (\kappa,r)\ \mapsto\ (\alpha,R)=(\arctan\kappa,\mathcal{R}(\kappa,r))
   \ \mapsto\ Re^{i\alpha},
\]
a composition of diffeomorphisms: $re^{i\beta}\mapsto(\tan\beta,r)$ is a diffeomorphism of $A$ onto
$\mathbb{R}\times\mathbb{R}_{>0}$; $(k,r)\mapsto(\kappa,r)$ is Lemma~\ref{isogone1} (real-analytic, with
$\partial_kT_1>0$); the map $(\kappa,r)\mapsto(\alpha,R)$ has Jacobian
$\det=\tfrac1{1+\kappa^2}\,\partial_r\mathcal {R}\ne0$,
since $\partial_r\mathcal {R}<0$ by Lemma~\ref{isogone2}(i), and its components are smooth because
$\mathcal {R}$ is (here $\mathcal{C}=\Re\mathcal{M}>0$ by Lemma~\ref{isogone2}(ii), so
$\mathcal{M}\neq0$ and $R=|\mathcal{M}|$ is real-analytic); and the polar map $(\alpha,R)\mapsto Re^{i\alpha}$
is a diffeomorphism where $R>0$. Hence $\mathcal{M}$ is a smooth bijection with nowhere-vanishing Jacobian, so a
diffeomorphism of $\mathcal{V}$ onto $\mathcal{W}$ by the inverse function theorem.
\proofend

\begin{remark}[continued fraction, Weyl disk, and boundary]\label{rem:weyl}
\emph{The bound as a continued fraction.} The disk bound of Lemma~\ref{isogone2}(ii) can be read
directly off the recurrence, which exposes why it is a half-plane condition. Writing $\nu=1-ik$, so that
$1/\mathcal{M}=\frac r2\,I_{\nu-1}(r)/I_\nu(r)$, the recurrence~\eqref{rec2} applied at
$\mu=\nu,\nu+1,\dots$ expands the ratio as the Gauss continued fraction
\[
   \frac1{\mathcal{M}}=a_0+\cfrac{x}{a_1+\cfrac{x}{a_2+\cfrac{x}{a_3+\cdots}}},\qquad
   a_n=(n+1)-ik,\quad x=(r/2)^2>0 .
\]
Its $N$-th approximant is $\phi_0\circ\cdots\circ\phi_{N-1}(a_N)$ with $\phi_n(t)=a_n+x/t$, and for
$\Re t>0$ one has $\Re\phi_n(t)=(n+1)+x\Re t/|t|^2>n+1\ge1$; thus each $\phi_n$ maps the right half-plane
into $\{\Re\ge1\}$, and composing inward from $a_N$ as $N\to\infty$ gives $\Re(1/\mathcal{M})\ge1$. The
bound is strict for $r>0$: the outermost step is $\Re(1/\mathcal{M})=1+x\Re(1/w_1)/|w_1|^2$ with tail
$w_1=\frac r2\,I_\nu/I_{\nu+1}$ of the same type, $\Re w_1\ge2$.

\emph{The Weyl reading, and a caveat.} The continued fraction converges to $1/\mathcal{M}$ by Pincherle's
theorem~\cite{Gautschi}: among the solutions of $\tfrac r2(c_{n+1}-c_{n-1})+(n-ik)c_n=0$ it selects the
ratio of the recessive one, here $c_n\propto I_{n-ik}(r)$, which by Proposition~\ref{prop:cont} are the
EMvM moments up to an $n$-independent factor, with $\mathcal{M}=\tfrac2r\,c_1/c_0$. This is the
Weyl--Titchmarsh $m$-function built from the recessive solution in the usual way~\cite{Teschl}, and the
half-plane $\{\Re\ge1\}$ is its value region, the analogue of a Weyl disk mapped to $W$ by $z\mapsto1/z$.
The caveat is that the transport term $\tfrac r2(c_{n+1}-c_{n-1})$ makes the off-diagonal \emph{skew}, so
the underlying Jacobi operator is not self-adjoint; accordingly $1/\mathcal{M}$ is not a
Herglotz--Nevanlinna function and the value region carries no spectral measure. What survives is the
geometry: the containment $\Re(1/\mathcal{M})\ge1$ above, and the radius $\tfrac12=\tfrac1{2\cdot1}$ set
by the order gap $1$ between $I_{1-ik}$ and $I_{-ik}$ (a $p$-th moment $I_{p-ik}/I_{-ik}$ would give
$\{\Re\ge p\}$ and the disk $\{|z-\tfrac1{2p}|\le\tfrac1{2p}\}$). This places the image \emph{inside} $W$;
the radial monotonicity of Lemma~\ref{isogone2}(i) is the extra ingredient fixing it as the \emph{entire}
disk.

\emph{The boundary.} The boundary circle $\partial W=\{R=\cos\alpha\}$ is the small--radius envelope of
$\mathcal{M}$: as $r\to0^+$, Proposition~\ref{prop:cont} and $I_\nu(r)\sim(r/2)^\nu/\Gamma(\nu+1)$ give
$\mathcal{M}(k,r)\to(1-ik)^{-1}$ (also the head $a_0$ of the continued fraction, at $x=0$), which runs
over $\partial W\setminus\{0\}$ with $\arg=\arctan k$ and modulus $(1+k^2)^{-1/2}=\cos(\arctan k)$. The
monotonicity of Lemma~\ref{isogone2}(i) then fills the disk inward from this circle to the star center
$0=\mathcal{M}(k,\infty)$, extending $\mathcal{M}$ continuously to $\partial W$ as $r\to0^+$. The von
Mises diameter is the segment $[0,1]$, traversed by the symmetric isogone $k=0$ as $r$ runs from $\infty$
to $0$.
\end{remark}

\proof[Proof of Theorem~\ref{main0}]
The mean-field self-consistency condition (\ref{OP3}) is always satisfied by the trivial solution $r=0$. Denoting $R=\frac2{\mu}$, the nontrivial solutions are determined from
\be\label{OP7}
\mathcal{M}(k,r)=Re^{i\alpha}, 
\ee
or, equivalently, from the system   
\begin{align}\label{OP4}
\mathcal{R}(k,r)&=R\\
\mathcal{A}(k,r)&=\alpha. \nonumber
\end{align}
By the above lemma, \eqref{OP7} has a unique solution as long as $Re^{i\alpha}\in W$. This means that for a fixed $\alpha\in(-\pi/2,\pi/2)$, there will exist a unique nontrivial solution $(k_\alpha(\mu),r_\alpha(\mu))$ 
of (\ref{OP3}) as long as $\mu>2\sec\alpha$. It is easy to see that the curve $\Gamma_\alpha$ consisting of nontrivial solutions from the statement of Theorem~\ref{main0} coincides exactly with $\Gamma_{\tan\alpha}$ from Lemma~\ref{isogone1}, while $K_\alpha(r)$ coincides exactly with $K(\tan\alpha,r)$. The fact that the bifurcation is global follows from the lack of horizontal or vertical asymptotes to $\Gamma_\alpha$, which is justified by the limits (\ref{lim1}) and (\ref{lim2}). The remaining statements follow immediately.  
\proofend  

\section*{Acknowledgments} 
The author would like to thank Professor Peter Constantin for the invaluable role he has played in his career, and for setting him on the path of exploring synchronization in nature.  

The author would like to thank Professor Andrew Poje for crucial insights and continuous support throughout this project.

\end{document}